\documentclass{amsart}
\usepackage{amsmath,amssymb,dsfont,amsthm,amsfonts}
\usepackage{graphics}
\usepackage{xypic}
\usepackage{latexsym}

\newtheorem{theo}{Theorem}

\renewcommand{\Im}{{\rm Im}}

\newcommand{\beq}{\begin{equation}}
\newcommand{\eeq}{\end{equation}}

\newcommand{\R}{{\mathds{R}}}

\newcommand{\A}{{\mathcal{A}}}
\renewcommand{\H}{{\mathcal{H}}}

\renewcommand{\b}{{\partial}}
\newcommand{\Vol}{{\rm Vol}}
\date{}

\title{Reidemeister torsion and analytic torsion of spheres}

\author[T. de Melo and M. Spreafico]{Thiago de Melo and Mauro Spreafico}

\address{\tt ICMC, Universidade de S\~ao Paulo, S\~ao Carlos, Brazil.}
\email{tmelo@icmc.usp.br, mauros@icmc.usp.br}



\begin{document}
\maketitle

\begin{abstract} We provide a simple topological derivation of a formula for the Reidemeister and the  analityc torsion of spheres.
\end{abstract}

\vskip .2in

Weng and You gave in \cite{WY} a  formula for the analytic torsion of a sphere with the standard metric. Their result is obtained by direct
calculation applying the definition of analytic torsion, and is based on the explicit knowledge of the spectrum of
the Laplacian on forms. Beside simplicity and generality, the result is not particularly illuminating, as stated in \cite{WY}. 
We prove in this note a
different formula for the analytic torsion of a sphere, and we show that this new formula is equivalent to the one given in \cite{WY}. Our proof is
based on purely geometric topological means, namely evaluation of the Reidemeister (R) torsion \cite{Rei} and application of the
Cheeger-M\"uller theorem \cite{Che} \cite{Mul}. The main motivations, beside the different proof, are: from one side, the topological approach is natural and simpler, from the other, our formula provides a nice geometric interpretation of the result, that now we state, after introducing some notation.  For a closed connected Riemannian manifold $(W,g)$, with metric $g$, and an orthogonal representation $\rho$ of the fundamental group of $W$, we denote by $T((W,g);\rho)$ the analytic torsion of $(W,g)$ with respect to $\rho$ (see \cite{RS} Definition 1.6), and by $\tau_{\rm R}((W,g);\rho)$ the Ray and   Singer version of the Reidemeister torsion of $(W,g)$ with respect to the same representation (see equation (\ref{e2}) below for the definition).

\begin{theo} \label{t1}
Let $S_l^n$ be the sphere of radius $l$ in $\R^{n+1}$ ($n>0$) with the standard Riemannian metric $g_l$ induced by the immersion. Let $\rho_0$ be a trivial orthogonal representation of $\pi_1(S^n_l)$. Then,
\[
T((S_l^n,g_l);\rho_0)=\tau_{\rm R}((S^n_l,g_l);\rho_0)=\left\{\begin{array}{cl}  1, & \textrm{\rm if $n$ is even},\\
( \Vol_{g_l}(S^n_l))^{{\rm rank}(\rho_0)}, & \text{\rm if $n$ is odd}.\end{array}\right.
\]
\end{theo}

The rest of this note is dedicated to the proof of this theorem. 
We recall that, by the Cheeger-M\"{u}ller theorem \cite{Che} \cite{Mul}, the analytic torsion $T((S_l^n,g_l);\rho)$ coincides with the R torsion $\tau_{\rm R}((S^n_l,g_l);\rho)$. Thus, it remains to evaluate the R torsion. For we first introduce some notation. Let
$$
\xymatrix{
C:& C_n\ar[r]^{\b_n}&C_{n-1}\ar[r]^{\b_{n-1}}&\dots\ar[r]^{\b_2}&C_1\ar[r]^{\b_1}&C_0,
}
$$
be a chain complex of real vector spaces. Denote by $Z_q=\ker \b_q$, by $B_q=\Im \b_{q+1}$, and  by
$H_q(C)=Z_q/B_q$ as usual.  For two bases $x=\{x_1,\dots, x_m\}$ and
$y=\{y_1,\dots, y_m\}$ of a vector space $V$, denote by $(y/x)$ the matrix defined by the change of base. For each $q$, fix a base $c_q$ for $C_q$, and a base $h_q$ for $H_q(C)$. Let $b_q$ be a set of (independent) elements of $C_q$ such that
$\b_q(b_q)$ is a base for $B_{q-1}$. Then the set of elements $\{\b_{q+1}(b_{q+1}),h_q, b_q\}$ is a base for $C_q$. In this situation, the Reidemeister torsion of the complex $C$ with respect to the graded base $h=\{h_q\}$ is the positive real number
\[
\tau_{\rm R}(C;h)=\prod_{q=0}^n \left| \det(\b_{q+1}(b_{q+1}),h_q,b_q/c_q)\right|^{(-1)^q}.
\]

Let $K$ be a connected finite cell complex of dimension $n$ and $\tilde K$ its universal covering complex, and identify the fundamental group of $K$
with the group of the covering transformations of $\tilde K$.  Let $C(\tilde K;\R)$ be the real chain complex of $\tilde K$. The action of the
group of covering transformations makes each chain group $C_q(\tilde K;\R)$ into a module over the group algebra $\R\pi_1(K)$, and each of these modules
is $\R \pi_1(K)$-free and finitely generated by the natural choice of the $q$-cells of $K$. 
We have got the complex $C(\tilde K;\R\pi_1(K))$ of free finitely generated modules over $\R \pi_1(K)$. Let $\rho:\pi_1(K)\to O(m,\R)$ be a representation of the fundamental group, and consider the twisted complex $C(K;\R^m_\rho)$. Assume $H_q(C(K; \R^m_\rho))$ are free finitely generated modules over $\R \pi_1(K)$. The Reidemeister torsion of $K$ with respect to the representation $\rho$ and to the graded base $h$ is defined applying the previous construction to the twisted complex $C(K;\R^m_\rho)$, namely 
\beq\label{e1}
\tau_{\rm R}(K;h,\rho)=\prod_{q=0}^n \left| \det\rho(\b_{q+1}(b_{q+1}),h_q,b_q/c_q)\right|^{(-1)^q},
\eeq
in $\R^+$. If $K$ is the cellular (or simplicial)
decomposition of a space $X$, the Reidemeister  torsion of $X$ is defined accordingly, and denoted by $\tau_{\rm R}(X;h,\rho)$. It was proved in
\cite{Mil} that $\tau_{\rm R}(X;h,\rho)$ does not depend on the decomposition $K$.

Let $W$ be a closed connected orientable Riemannian manifold of dimension $n$ with Riemannian metric $g$. Then, all the previous assumptions are satisfied, and the R torsion $\tau_{\rm R}(W;h,\rho)$ is well defined for each fixed graded base $h$ for the homology of $W$, and each representation $\rho$ of the fundamental group. In this context, Ray and Singer suggest a natural geometric invariant object, by fixing an appropriate base $h$ using the geometric
structure, as follows.  Let $E_\rho\to W$ be the real vector bundle associated to the representation $\rho:\pi_1(W)\to O(m,\R)$. Let
$\Omega(W,E_\rho)$ be the graded linear space of smooth forms on $W$ with values in $E_\rho$.  The exterior differential on $W$ defines the exterior differential on $\Omega^q(W, E_\rho)$, $d:\Omega^q(W, E_\rho)\to
\Omega^{q+1}(W, E_\rho)$. The metric $g$ defines an Hodge operator on $W$ and hence on $\Omega^q(W, E_\rho)$, $\star:\Omega^q(W, E_\rho)\to
\Omega^{n-q}(W, E_\rho)$, and,  using the inner product in $E_\rho$, an inner product on $\Omega^q(W, E_\rho)$. Let $\H^q$ be the space of the harmonic
forms in $\Omega^q(W, E_\rho)$, and let $\A^q$ be the de Rham map $\A^q:\H^q\to C^q(W;E_\rho)$,
\[
\A^q(\omega)(c\otimes_\rho v)=\int_{ c} (\omega, v),
\]
with $c\otimes_\rho v\in C_q(W;E_\rho)$, and where we identify the chain $c$ with the $q$-cell that $c$ represents, and $(~,~)$ is the inner product in $\R^m$.
Following Ray and Singer, we define the following map, where $\hat c_{q,j}$ is the Poincar\'e dual cell to $c_{q,j}$,
\beq\label{aa}
\begin{aligned}
\A_q:&\H^q\to C_q(W;E_\rho),\\
\A_q:&\omega\mapsto(-1)^{(n-1)q}\sum_{j,k} \int_{\hat c_{q,j}}(\star\omega,e_k)
(c_{q,j}\otimes_\rho e_k).
\end{aligned}
\eeq

In this situation, let $a$ be a graded orthonormal base for the space of the harmonic forms in $\Omega (W,E_\rho)$, then we call the
positive real number
\beq\label{e2}
\tau_{\rm R}((W,g);\rho)=\tau_{\rm R}(W;\A(a),\rho)=\prod_{q=0}^n \left| \det\rho(\b_{q+1}(b_{q+1}),\A_q(a_q),b_q/c_q)\right|^{(-1)^q},
\eeq
the R torsion of $(W,g)$ with respect to the representation $\rho$. It can be proved that  $\tau_{\rm R}((W,g);\rho)$ does not depend
on the choice of the orthonormal base $a$.

We now compute the R torsion of the sphere $S^n_l$. Recall that the sphere $S^n_l=\{x\in \R^{n+1}~|~|x|=l\}$, is parameterized in polar coordinates
by
$$\left\{
\begin{array}{rcl}
x_1&=&l\sin{\theta_n}\sin{\theta_{n-1}}\cdots\sin{\theta_3}\sin{\theta_2}\sin{\theta_1} \\[8pt]
x_2&=&l\sin{\theta_n}\sin{\theta_{n-1}}\cdots\sin{\theta_3}\sin{\theta_2}\cos{\theta_1} \\[8pt]
x_3&=&l\sin{\theta_n}\sin{\theta_{n-1}}\cdots\sin{\theta_3}\cos{\theta_2} \\[8pt]
x_4&=&l\sin{\theta_n}\sin{\theta_{n-1}}\cdots\cos{\theta_3} \\[8pt]
&\vdots& \\
x_n&=&l\sin{\theta_n}\cos{\theta_{n-1}} \\[8pt]
x_{n+1}&=&l\cos{\theta_n}
\end{array}
\right.$$
with
$0\leq\theta_1\leq 2\pi$, $0\leq \theta_2,\ldots,\theta_n\leq \pi$. The induced Riemannian metric is
\[
g_l=l^2\left((d\theta_n)^2+\sin^2\theta_n(d\theta_{n-1})^2+\cdots+\prod_{i=2}^n\sin^2\theta_i(d\theta_1)^2\right)=l^2 g_1.
\]

Let $K$ be the standard cellular decomposition of $S^n_l$, with one top cell and one $0$-cell. Let $\rho_0$ be the trivial representation when $n=1$, of rank $m$. Then the relevant complex is
\[
C: \quad \xymatrix{0\ar[r] & \R^m[c_n^1] \ar[r] & 0\ar[r]&\cdots \ar[r]&0\ar[r] & \R^m[c_0^1] \ar[r] & 0},
\]
with preferred bases $c_0=\{c_0^1\}$ and  $c_n=\{c_n^1\}$.  To fix the base for the homology, we need a graded  orthonormal base $a$ for the
harmonic forms. This is $a_0=\left\{\frac{1}{\sqrt{\Vol_{g_l}(S^n_l)}}\right\}$ and ${a}_n=\left\{\frac{1}{\sqrt{\Vol_{g_l}(S^n_l)}}
\sqrt{|g_l|}d\theta_n\wedge\cdots\wedge d\theta_1\right\}$. Applying the formula in equation (\ref{aa}) we obtain $h_0=\{h_0^1\}$, $h_n=\{h_n^1\}$,
with
\begin{align*}
h^1_0 &=\A_0(a^1_0)=\frac{1}{\sqrt{\Vol_{g_l}(S^n_l)}}\int_{S^n_l} \sqrt{|g_l|}d\theta_n\wedge\cdots\wedge
 d\theta_1 c_0^1= \sqrt{\Vol_{g_l}(S^n_l)} c_0^1, \\
h^1_n &=\A_n(a^1_n)= \frac{1}{\sqrt{\Vol_{g_l}(S^n_l)}}\int_{\rm pt}\star \sqrt{|g_l|}d\theta_n\wedge\cdots\wedge
 d\theta_1 c^1_n= \frac{1}{\sqrt{\Vol_{g_l}(S^n_l)}}  c_n^1.
 \end{align*}

As $b_q=\emptyset$, for all $q$, we have that
\[
|\det\rho (h_0/c_0)| = \left(\sqrt{\Vol(S^n_l)}\right)^m, \qquad |\det \rho(h_n/c_n)| = \frac{1}{\left(\sqrt{\Vol(S^n_l)}\right)^m}.
\]

Applying the definition in equation (\ref{e2}), this gives
\begin{align*}
\tau_{\rm R}((S^n_l,g_l);\rho_0)  = & \left(\sqrt{\Vol_{g_l}(S^n_l)}\right)^m \left(\frac{1}{\sqrt{\Vol_{g_l}(S^n_l)}}\right)^{m(-1)^n}, \\
\end{align*}
and  completes  the proof of the theorem.

We conclude with some remarks.
\begin{itemize}

\item The formula in Theorem \ref{t1} is an expected result for spheres. In fact, it is well known that fixing non zero volume elements $u_q\in
\Lambda^{d_q}(H^q(W;E_\rho))$, $d_q=\dim H^q(W;E_\rho)$, and hence fixing a graded base, say $h(u)$ for the homology, we have
\[
\tau_{\rm R}(W;h(\alpha u))=\frac{\alpha_0}{\alpha_1}\frac{\alpha_2}{\alpha_3}\cdots \tau_{\rm R}(W;h(u)),
\]
where $\alpha u=\{\alpha_q u_q\}$, $\alpha_q\in \R$.

\item If we use the cellular decomposition with two cells in each dimension, a similar calculation gives the same result, since the determinant in
all dimensions $q$, with $0<q<n$, is $1$.

\item Consider the product $S^n_a\times S^k_b$ with the product metric. Using the cellular decomposition with four cells: $c_0^1$, $c_n^1$, $c_k^1$,
and $c_{n+k}^1$, we obtain
\[
T(S^n_a\times S^k_b)=\left\{\begin{array}{cl}
    \Vol(S^k_{b})^{\chi(S^n_{a})} & \quad \text{$n$ even, $k$ odd}, \\[10pt]
    \Vol(S^n_{a})^{\chi(S^k_{b})} & \quad \text{$k$ even, $n$ odd}, \\[10pt]
    1 & \quad \text{$n,k$ even or $n,k$ odd},\end{array}\right.
\]
in agreement with one of the main properties of the torsion.

\item We prove the equivalence of the formula given in Theorem \ref{t1} (in the case of ${\rm rank}(\rho_0)=1$) with the formula given in \cite{WY}. This follows recalling the volume of the 
sphere $S_l^n$
\[
\Vol_{g_l} (S^n_l)=\frac{2\pi^\frac{n+1}{2}l^n}{\Gamma\left(\frac{n+1}{2}\right)},
\]
where $\Gamma(z)$ is the Euler Gamma function. This gives:
\[
T((S_l^{2k+1},g_l);\rho_0)=\frac{2\pi^{k+1}l^{2k+1}}{k!}.
\]

\end{itemize}

{\bf Acknowledgments.} M. Spreafico thanks M. Lesh for useful conversations.

\end{document}